\documentclass[12pt]{amsart}
\usepackage{amssymb}
\usepackage{amsbsy}
\usepackage{amscd}

\usepackage[mathscr]{eucal}
\usepackage{mathbbol}

\DeclareSymbolFontAlphabet{\mathbb}{AMSb}
\DeclareSymbolFontAlphabet{\mathbbl}{bbold}

\DeclareSymbolFont{rsfs}{U}{rsfs}{m}{n}
\DeclareSymbolFontAlphabet{\mathrsfs}{rsfs}

\usepackage{a4wide}
\usepackage[all]{xy}
\usepackage{tikz-cd}
\usepackage{tikz}
\usetikzlibrary{arrows,decorations.pathmorphing,backgrounds,positioning,fit,petri}
\usetikzlibrary{decorations.markings}
\usetikzlibrary{patterns,calc}
\usepackage{geometry}
\usetikzlibrary{matrix,positioning,decorations.pathreplacing}

\usepackage{verbatim}
\usepackage{version}
\usepackage{color}
\definecolor{darkspringgreen}{rgb}{0.09, 0.45, 0.27}
\definecolor{deepjunglegreen}{rgb}{0.0, 0.29, 0.29}
\newenvironment{NB}{
\color{red}{\bf NB}. \footnotesize
}{}

\excludeversion{NB}
\excludeversion{NB2}
\excludeversion{NB3}
\makeatletter

\usepackage{url}
\usepackage{ifmtarg}

\usepackage{xr-hyper}
\usepackage[
pdftex,
bookmarks=true,
colorlinks=true,
citecolor=deepjunglegreen,
filecolor=deepjunglegreen,
debug=true,
unicode,
pdfnewwindow=true]{hyperref}
\usepackage{cleveref}

\crefname{Theorem}{Theorem}{Theorems}
\crefformat{Theorem}{Theorem~#2#1#3}
\crefname{section}{\S}{\S\S}
\crefformat{section}{\S#2#1#3} 
\crefmultiformat{section}{\S\S#2#1#3}{, #2#1#3}{, #2#1#3}{, #2#1#3}
\crefname{Lemma}{Lemma}{Lemmas}
\crefformat{Lemma}{Lemma~#2#1#3}
\crefname{Proposition}{Proposition}{Propositions}
\crefformat{Proposition}{Proposition~#2#1#3}
\crefname{Corollary}{Corollary}{Corollaries}
\crefformat{Corollary}{Corollary~#2#1#3}
\crefname{Definition}{Definition}{Definitions}
\crefformat{Definition}{Definition~#2#1#3}
\crefname{Remark}{Remark}{Remarks}
\crefformat{Remark}{Remark~#2#1#3}
\crefname{Remarks}{Remark}{Remarks}
\crefformat{Remarks}{Remark~#2#1#3}
\crefname{Conjecture}{Conjecture}{Conjectures}
\crefformat{Conjecture}{Conjecture~#2#1#3}
\crefname{figure}{Figure}{Figure}
\crefformat{figure}{Figure~#2#1#3}

\crefname{appendix}{Appendix}{Appendices}
\crefformat{appendix}{\S#2#1#3}
\crefformat{enumi}{(#2#1#3)}

\crefname{equation}{}{}

\usepackage{enumitem}
\usepackage{youngtab}
\usepackage{mathtools}
\usepackage{extpfeil}
\usepackage{blkarray}

\usepackage{multirow,bigdelim}
\usepackage{array}

\renewcommand{\thesubsection}{\thesection(\@roman\c@subsection)}
\newcounter{number}
\makeatother
\newtheorem{Theorem}[equation]{Theorem}

\theoremstyle{definition}

\theoremstyle{remark}

\numberwithin{equation}{section}

\newcommand{\defeq}{\overset{\operatorname{\scriptstyle def.}}{=}}
\newcommand{\CC}{{\mathbb C}}
\newcommand{\ZZ}{{\mathbb Z}}

\newcommand{\RR}{{\mathbb R}}

\newcommand{\SL}{\operatorname{\rm SL}}

\newcommand{\GL}{\operatorname{GL}}

\newcommand{\algsl}{\operatorname{\mathfrak{sl}}} 

\newcommand{\Spec}{\operatorname{Spec}\nolimits}

\newcommand{\Hom}{\operatorname{Hom}}

\renewcommand{\MR}[1]{}

\newcommand{\vin}[1]{\operatorname{i}(#1)} 
\newcommand{\vout}[1]{\operatorname{o}(#1)} 

\newcommand{\bN}{\mathbf N}

\newcommand{\shfO}{\mathcal O}

\newcommand{\tslabar}{\mathbin{
\setbox0=\hbox{/\!\!/\!\!/}\rule[0.4\ht0]{\wd0}{.3\dp0}\kern-\wd0\box0}}

\newcommand{\aff}{\mathrm{aff}}

\newcommand{\Gr}{\mathrm{Gr}}

\newcommand{\cR}{\mathcal R}

\newcommand{\cK}{\mathcal K}
\newcommand{\cO}{\mathcal O}
\newcommand{\cW}{\mathcal W}

\makeatletter
\newcommand{\cA}[1][{}]{%
  \@ifmtarg{#1}%
  {\mathcal A}
  {\mathcal A(#1)}
}
\newcommand{\cAh}[1][{}]{%
  \@ifmtarg{#1}%
  {\mathcal A_\hbar}
  {\mathcal A_\hbar(#1)}
}
\makeatother

\usepackage{accents}
\makeatletter
\newcommand{\sbullet}{%
  \hbox{\fontfamily{lmr}\fontsize{.4\dimexpr(\f@size pt)}{0}\selectfont\textbullet}}

\makeatother

\newcommand{\ft}{\mathfrak t}

\newcommand{\po}{\ar@{}[dr]|{\text{\pigpenfont R}}}
\newcommand{\pb}{\ar@{}[dr]|{\text{\pigpenfont J}}}
\newcommand{\pp}{\ar@{}[dr]|{\text{\pigpenfont P}}}

\newcommand{\cM}{\mathcal M}

\newcommand{\cS}{{\mathcal{S}}}

\newcommand{\oW}{\overline{\mathcal{W}}{}}

\newcommand{\fg}{{\mathfrak{g}}}

\newcommand{\fu}{{\mathfrak{u}}}

\makeatletter
\newdimen\y@inside
\y@inside=\y@boxdim
\advance\y@inside by -.2\y@linethick
\definecolor{halfgray}{gray}{0.75}
\newcommand{\rf}{\color{halfgray}
\rule[-.2\y@inside]{\y@inside}{\y@inside}
}
\makeatother

\DeclareSymbolFont{symbolsC}{U}{pxsyc}{m}{n}
\SetSymbolFont{symbolsC}{bold}{U}{pxsyc}{bx}{n}
\DeclareFontSubstitution{U}{pxsyc}{m}{n}
\DeclareMathSymbol{\medcirc}{\mathbin}{symbolsC}{7}

\newcommand{\fin}{{\operatorname{fin}}}

\newcommand{\IC}{\operatorname{IC}}

\newcommand{\relmiddle}[1]{\mathrel{}\middle#1\mathrel{}}

\newcommand{\Heis}{\operatorname{Heis}}

\begin{document}

\title[Instanton moduli spaces and affine flag varieties]
{Intersection cohomology groups of instanton moduli spaces
  and cotangent bundles of affine flag varieties}
\author[H.~Nakajima]{Hiraku Nakajima}
\address{Kavli Institute for the Physics and Mathematics of the Universe (WPI),
  The University of Tokyo,
  5-1-5 Kashiwanoha, Kashiwa, Chiba, 277-8583,
  Japan
}
\email{hiraku.nakajima@ipmu.jp}
\maketitle
\pagestyle{empty}
\thispagestyle{empty}

This is an abstract for my talk at the 68th Geometry Sympoisum on
August 31, 2021. It is based on my joint work in progress with {\bf
  Dinakar Muthiah}.

\section{Affine Grassmannian and the cotangent bundle of the flag variety}

Let us first recall a result of Arkhipov-Bezrukavnikov-Ginzburg
\cite{MR2053952}.
Let $G$ be a simply-connected complex semisimple group with the Lie
algebra $\fg$. We fix a Borel subgroup $B$ of $G$ and the maximal
torus $T\subset B$. Let $G/B$ be the flag variety and $T^*(G/B)$ its
cotangent bundle. Let $\mathfrak b$ be the Lie algebra of $B$. Then
$T^*(G/B) = G\times_B (\fg/\mathfrak b)^*$.
Let $\mu$ be a dominant weight of $G$ and consider a line bundle
$\shfO(\mu)$ over $T^*(G/B)$, pull-backed from $G/B$. Let us
consider the space
\(
  H^0(T^*(G/B), \shfO(\mu))
\)
of sections of $\shfO(\mu)$.

On the other hand, let us consider the affine Grassmannian
$\Gr_{G^\vee}$ associated with the Langlands dual group $G^\vee$,
which is an infinite dimensional partial flag variety $G^\vee(\cK)/G^\vee(\cO)$.
Here $\cO = \CC[[z]]$ (formal power series ring) and
$\cK = \CC((z))$ (formal Laurent power series ring).
Let us denote by $T^\vee$ a maximal torus of $G^\vee$, which is the
dual torus of $T$.
We consider a coweight $\lambda$ as an element of $T^\vee(\cK)$, and
further a point $z^\lambda\in G^\vee(\cK)/G^\vee(\cO)$ by the
composite $T^\vee(\cK)\to G^\vee(\cK)\to \Gr_{G^\vee}$.
Note that $\lambda$ can be also regarded as a weight of $G$.
We have an action of $G^\vee(\cO)$ on $\Gr_{G^\vee}$. It is well-known
that $G^\vee(\cO)$-orbits are parametrized by dominant coweights
$\lambda$ of $G^\vee$ as $G^\vee(\cO)z^\lambda$. Let
$\overline{\Gr}_{G^\vee}^\lambda$ denote the closure of
$G^\vee(\cO)z^\lambda$. It is known that
$\overline{\Gr}_{G^\vee}^\lambda$ is a projective variety. It is an
analog of a Schubert variety in the (finite dimensional) flag
variety. We consider the intersection cohomology complex
$\IC(\overline{\Gr}_{G^\vee}^\lambda)$ associated with the constant
local system.

We consider the above $\mu$ as a dominant coweight of $G^\vee$, and
denote the inclusion of the point $z^\mu$ in $\Gr_{G^\vee}$ by
$i_\mu$. Then one of main results in \cite{MR2053952} is a
$G$-equivariant graded vector space isomorphism
\begin{equation}\label{eq:1}
  H^0(T^*(G/B), \shfO(\mu)) \cong \bigoplus_\lambda V(\lambda)^*
  \otimes H^*(i_\mu^! \IC(\overline{\Gr}_{G^\vee}^\lambda))[\mu(2\rho^\vee)],
\end{equation}
where the summation runs over dominant weights $\lambda$.
We denote the irreducible representation of $G$ with the highest
weight $\lambda$ by $V(\lambda)$,
a cohomological degree shift by $[\ \ ]$,
\begin{NB}
  Contrary to the convention in \cite{MR3352528}.
\end{NB}%
and the sum of fundamental coweights of $G$ by $\rho^\vee$. The grading on
the left hand side is given by the $\CC^\times$-action, the dilatation
on fibers of $T^*(G/B)$. Moreover the isomorphism \eqref{eq:1} is
compatible with products: the product in the left hand side is given
by $\shfO(\mu_1)\otimes\shfO(\mu_2) \cong \shfO(\mu_1+\mu_2)$. The
product on the right hand side is given by the convolution on
$\Gr_{G^\vee}$.
See \cite[Th.~8.5.2]{MR2053952}.

We view \eqref{eq:1} as a construction of $T^*(G/B)$ by
\emph{topology} of the affine Grassmannian, as $T^*(G/B)$ can be
recovered from $\bigoplus_\mu H^0(T^*(G/B),\shfO(\mu))$ as a
homogeneous spectrum.
We will see similar constructions below.

For a later purpose let us extend \eqref{eq:1}. The maximal torus
$T^\vee$ acts on $\Gr_{G^\vee}$ so that $z^\mu$ is a fixed point.
Therefore we can consider the $T^\vee$-equivariant cohomology
$H^*_{T^\vee}(i_\mu^!\IC(\overline{\Gr}^\lambda_{G^\vee}))$. On the
other hand we replace $T^*(G/B) \cong G\times_B(\fg/\mathfrak b)^*$ by
a larger space $G\times_B(\fg/\mathfrak u)^*$ where $\mathfrak u$ is
the nilpotent radical of $\mathfrak b$. It gives Grothendieck's
simultaneous resolution $G\times_B(\fg/\mathfrak u)^*\to\fg^*$ of the
deformation $\fg^*$ of the nilpotent cone $\mathcal N$, parametrized
by $\mathfrak t^*$. Here $\mathfrak t^*$ is the dual of the Lie
algebra of $T$, which is the Lie algebra of $T^\vee$, and
$\Spec H^*_{T^\vee}(\mathrm{pt})$. We have
\begin{equation}\label{eq:2}
  H^0(G\times_B(\fg/\mathfrak u)^*, \shfO(\mu)) \cong
  \bigoplus_\lambda V(\lambda)^*
  \otimes H^*_{T^\vee}
  (i_\mu^! \IC(\overline{\Gr}_{G^\vee}^\lambda))[\mu(2\rho^\vee)].
\end{equation}
This is a $G$-equivariant isomorphism of graded
$\CC[\mathfrak t^*]\cong H^*_{T^\vee}(\mathrm{pt})$-modules.
It is specialized to \eqref{eq:1} at the fiber of $0$ in $\ft^*$.

\section{Singular monopoles}\label{sec:monopole}

Let us rewrite $H^*(i_\mu^! \IC(\overline{\Gr}_{G^\vee}^\lambda))$ so
that it looks a slightly familiar to geometers and physicists.

Let us first define a slice to the $G^\vee(\cO)$-orbit through
$z^\mu$. Let $G^\vee[z^{-1}]_1$ be the kernel of the evaluation at
$\infty$ homomorphism $G^\vee[z^{-1}]\to G^\vee$. Let
$\cW_{G^\vee,\mu}$ denote the orbit $G^\vee[z^{-1}]_1 z^\mu$ and
define
\begin{equation*}
  \oW^{\lambda}_{G^\vee,\mu} \defeq \overline{\Gr}_{G^\vee}^\lambda\cap\cW_{G^\vee,\mu}.
\end{equation*}
It is known that $\oW^\lambda_{G^\vee,\mu}$ is a transversal slice to
$\Gr^\mu_{G^\vee}$ in $\overline{\Gr}_{G^\vee}^\lambda$. Then we can replace 
$H^*(i_\mu^! \IC(\overline{\Gr}_{G^\vee}^\lambda))[\mu(2\rho^\vee)]$ by
\begin{equation*}
  H^*(i_\mu^! \IC(\oW^{\lambda}_{G^\vee,\mu})),
\end{equation*}
where $i_\mu^!$ is now regarded as an inclusion to
$\oW^{\lambda}_{G^\vee,\mu}$. Thus \eqref{eq:1} means that $T^*(G/B)$
is constructed by topology of $\oW^{\lambda}_{G^\vee,\mu}$.

The space, called \emph{affine Grassmannian slice} is also a
differential geometric object.
Let $G^\vee_c$ be the maximal compact subgroup of $G^\vee$. We
consider a $G^\vee_c$-monopole on $\RR^3$, i.e., a solution of the
Bogomolny equation for a pair $(A,\Phi)$, a $G^\vee_c$-connection and
a section of the adjoint bundle:
\begin{equation*}
  *F_A = d_A\Phi.
\end{equation*}
They satisfy the asymptototic condition at $|x|\to\infty$, given by
$\mu$. To make sense of this condition, $\mu$ must be a coweight of
$G^\vee$, but not necessarily dominant.
(We are considering the so-called maximal symmetry breaking
case.)

We further assume that $(A,\Phi)$ has a Dirac-type singularity at $0$,
which means that it comes from an $S^1$-equivariant instantons on
$\RR^4$ locally around $0$. Then the fiber over the origin in $\RR^4$
is equipped with a homomorphism $\tilde\lambda\colon S^1\to G^\vee_c$,
whose conjugacy class is the dominant coweight $\lambda$ of $G^\vee$
above.

Then we consider the framed moduli space $\cM^\lambda_{G^\vee,\mu}$ of
such singular $G^\vee_c$-monopoles. As for moduli spaces of instantons
on $4$-manifolds, we need to consider Uhlenbeck (partial)
compactification to get a reasonable space. Since we are considering
$S^1$-equivariant instantons, the bubbling could occur the origin in
$\RR^4$. Therefore the symmmetric product of the base $4$-manifold
does not appear, and we consider
\begin{equation*}
  \bigcup_{\mu\le\nu\le\lambda} \cM^\nu_{G^\vee,\mu},
\end{equation*}
where $\nu$ runs over the set of dominant coweights between $\mu$ and
$\lambda$.
Then this space is homeomorphic to $\oW^{\lambda}_{G^\vee,\mu}$. This
was proved by Braverman-Finkelberg \cite{braverman-2007}.
A very similar result, namely for local behaviour around $0$, was
explained earlier by Kapustin-Witten \cite{MR2306566}.

A little bit more precisely, Braverman-Finkelberg \cite{braverman-2007}
imposed a different boundary condition at $|x|\to \infty$ from the
usual one imposed for monopoles on $\RR^3$ such as in
\cite{MR709461,MR994495,MR934202}.
Considering a singular monopole on $\RR^3$ as an $S^1$-equivariant
instanton on $\RR^4$ as above, they imposed that the instanton extends
to $S^4 = \RR^4\cup\{\infty\}$.
Then $\mu$ is regarded as the homomorphism corresponding to the fiber
at $\infty$. Such $\mu$ is dominant.
Moreover it is expected that two monopole moduli spaces with different
asymptotic behavior are isomorphic as algebraic varieties (but
different hyper-K\"ahler structures) when $\mu$ is dominant and the
latter space can be defined.
We will return back this issue when we discuss Coulomb branches of
quiver gauge theories in \cref{sec:KacMoody}.

Let me also mention that similar results in closely related situations
were studied in \cite{MR2755487,MR3536929}.



\section{Geometric Satake correspondence}

The isomorphism \eqref{eq:1} was built on the geometric Satake
correspondence
\cite{Lus-ast,1995alg.geom.11007G,MV2,Beilinson-Drinfeld}, which gives
an equivalence of tensor categories
\begin{equation*}
  \operatorname{Rep}G \cong \operatorname{Perv}_{G^\vee(\cO)}\Gr_{G^\vee}.
\end{equation*}
Here $\operatorname{Rep}G$ is the category of finite dimensional
representations of $G^\vee$, which is a semisimple abelian category
equipped with the tensor product $\otimes$. The right hand side is the
category of $G^\vee(\cO)$-equivariant perverse sheaves on
$\Gr_{G^\vee}$, which is an abelian category equipped with the tensor
product givey by the convolution. Moreover this equivalence sends
$V(\lambda)$ to $\IC(\overline{\Gr}_{G^\vee}^\lambda)$.

The weight space $V(\lambda)_\mu$ of $V(\lambda)$ is realized as
follows (see \cite{MV2} and also \cite{braverman-2007} for an explicit
use of $\oW^\lambda_{G^\vee,\mu}$): Let $\rho$ be the sum of
fundamental weights of $G$. We consider it as a regular dominant
cocharacter $\rho\colon\CC^\times\to T^\vee$. We then consider the
induced $\CC^\times$-action on $\oW^\lambda_{G^\vee,\mu}$.
It is well-known that the fixed point set, denoted by
$(\oW^\lambda_{G^\vee,\mu})^{\rho=0}$, consists of a single point
$z^\mu$.
More generally we consider an affine algebraic variety $X$ with
a $\CC^\times$-action.
We define the \emph{repelling set}
\begin{equation*}
  X^{\rho\le 0}
  \defeq \left\{ x \in X \relmiddle|
  \lim_{t \rightarrow \infty} \rho(t)\cdot x \text{ exists} \right\}.
\end{equation*}
We have a diagram
\begin{equation*}
  \begin{tikzcd}
    X^{\rho=0} \arrow[shift left=.5,r,"\iota"] & X^{\rho\le 0} \arrow[r,"j"]
    \arrow[shift left=.5, l, "p"] & X,
  \end{tikzcd}
\end{equation*}
where $\iota$, $j$ are the inclusion, and $p$ denotes the morphism
given by taking the limit $\eta(t)\cdot x$ for $t\to \infty$.
Then the hyperbolic restriction functor $\Phi$ is defined as
$\iota^* j^!$ on objects on the $\CC^\times$-equivariant derived
category. See \cite{Braden,MR3200429}.
Then we have
\begin{equation*}
  V(\lambda)_\mu \cong H^0(\Phi(\IC(\oW^\lambda_{G^\vee,\mu}))).
\end{equation*}
In this case it is also known that the hyperbolic restriction functor
vanishes in other degrees than $0$.

Note that there is a difference between this cohomology and one
appears in \eqref{eq:1}, i.e.,
$H^*(i_\mu^! \IC(\overline{\Gr}_{G^\vee}^\lambda))[\mu(2\rho^\vee)]
\cong H^*(i_\mu^! \IC(\oW_{G^\vee,\mu}^\lambda))$. It is the costalk
at $z^\mu$. On the other hand, the above is defined by
$\Phi = \iota^* j^!$. Compare it with $i_\mu^! = \iota^! j^!$.
They are same dimensional vector spaces, it is hard to distinguish
them.
We will clarify the difference in \cref{sec:main}. In the situation of
the geometric Satake for the finite dimensional group $G$, it was
clarified by Ginzburg-Riche \cite{MR3352528}. We will follow it
in the Kac-Moody setting.

\section{Coulomb branches and ring objects}\label{sec:Coulomb}

Let us return to \eqref{eq:1}. We rewrite the right hand side as
\begin{equation*}
  H^*(i_\mu^! \mathcal A_R)[\mu(2\rho^\vee)]\quad
  \text{where } \mathcal A_R = \bigoplus_\lambda V(\lambda)^*\otimes
  \IC(\overline{\Gr}_{G^\vee}^\lambda).
\end{equation*}
Then observe that $\mathcal A_R$ corresponds to
$\bigoplus_\lambda V(\lambda)^*\otimes V(\lambda)$ under the geometric
Satake correspondence. It is further isomorphic to the regular
representation $\CC[G]$ by the Peter-Weyl theorem. Note $\CC[G]$ is a
commutative ring object in $\operatorname{Rep}(G)$ (more precisely its
completion as we take an infinite sum). The concept of a ring object
is an abstract notion: it is equipped with a homomorphism
$\CC[G]\otimes\CC[G]\to \CC[G]$ satisfying associativity,
commutativity, etc. Therefore $\mathcal A_R$ is a ring object in
$\operatorname{Perv}_{G^\vee(\cO)}\Gr_{G^\vee}$.

This observation gave motivation to study more general \emph{ring
  objects} in $\operatorname{Perv}_{G^\vee(\cO)}\Gr_{G^\vee}$, or the
corresponding derived category $D_{G^\vee(\cO)}\Gr_{G^\vee}$ in
\cite{2017arXiv170602112B}. Given a ring object $\mathcal A$ in
$D_{G^\vee(\cO)}\Gr_{G^\vee}$, we apply the same construction to get a
graded comutative ring. We further take its spectrum or homogeneous
spectrum to define a scheme. The result could be a new scheme. Even if
we get a known scheme, e.g. $T^*(G/B)$ as above, we could get a new
perspective to the known scheme.

The mathematical definition of the Coulomb branch of a $3D$ SUSY gauge
theory \cite{2015arXiv150303676N,2016arXiv160103586B} can be regarded
as a construction of a ring object. Let us briefly recall the
construction. Suppose a representation $\bN$ of $G^\vee$ is given. We
consider the \emph{variety $\mathcal R$ of BFN triples}
$(\mathcal P,\varphi,s)$:
\begin{description}
\item[B{\rm undle}] a principal $G^\vee$-bundle $\mathcal P$ over the formal disk
  $D = \Spec\cO$,
\item[F{\rm raming}] a framing
  $\varphi\colon \mathcal P|_{D^*} \to G\times D^*$ of $\mathcal P$
  over the formal punctured disk $D^* = \Spec\cK$,
\item[N{\rm-valued section}] a section $s$ of an associated vector
  bundle $\mathcal P\times_{G^\vee}\bN$ over $D$.
\end{description}
We impose the condition that $\varphi(s)$ extends across $0\in D$. It
is a nontrivial constraint on $s$ as $\varphi$ is defined only over $D^*$.

If we only consider a pair $(\mathcal P,\varphi)$, it defines the
affine Grassmannian $\Gr_{G^\vee}$.

A convolution product is defined on the equivariant Borel-Moore
homology $H^{G^\vee(\cO)}_*(\cR)$, which is shown to be commutative. Then
the Coulomb branch for a SUSY gauge theory associated with
$(G^\vee,\bN)$ is defined as $\Spec H^{G^\vee(\cO)}_*(\cR)$.

Let $\omega_\cR$ be the dualizing complex on $\cR$. Then
$H^{G^\vee(\cO)}_*(\cR) \cong H^*_{G^\vee(\cO)}(\omega_\cR)$. The cohomology
$H^*_{G^\vee(\cO)}$ can be regarded as pushforward for the projection
$\cR\to\mathrm{pt}$. We forget $\bN$-sections from the variety of
triples, we have the morphism $\pi\colon \cR\to\Gr_{G^\vee}$.
The convolution product on $H^{G^\vee(\cO)}_*(\cR)$ also gives a structure
of a ring object on $\pi_*\omega_\cR$. See \cite{2017arXiv170602112B}.

We explain two nontrivial sources of ring objects in
$D_{G^\vee}(\Gr_{G^\vee})$ so far. We further have several
manupulations constructing new ring objects from given ones. Therefore
we get many ring objects, which give interesting spaces, such as Higgs
branches of class $S$ theories asked by Moore-Tachikawa
\cite{MR2985331}.

\section{Geometric Satake correspondence for a Kac-Moody Lie
  algebra}\label{sec:KacMoody}


Braverman-Finkelberg \cite{braverman-2007} proposed analog of
geometric Satake correspondence for an affine Kac-Moody group $G$
(assumed to be dual of untwisted and simply-connected type).
They used Uhlenbeck partial
compactification of $(G^\vee)_{\fin,c}$-instantons on
$\cS_\ell \defeq \RR^4/(\ZZ/\ell\ZZ)$, where $G^\vee$ is the dual of
$G$, $(G^\vee)_\fin$ is the associated finite dimensional group,
$(G^\vee)_{\fin,c}$ is a maxmimal compact subgroup of $(G^\vee)_\fin$,
and $\ell$ is the level of representations of $G$.
Uhlenbeck space serves as an analog of $\oW^\lambda_{G^\vee,\mu}$. It
is a natural analog as monopoles are $S^1$-equivariant instantons as
we explained above.

When $G$ is of type $A^{(1)}_{n-1}$, proposed statements were mostly
checked by the second-named author \cite{Na-branching}. Uhlenbeck
spaces are quiver varieties of type $A_{\ell-1}^{(1)}$, hence related
to representations of $\algsl(\ell)_\aff$
\cite{Na-quiver,Na-alg}. Then level rank duality
\cite{MR675108,Hasegawa} between $\algsl(n)_\aff$ at level $\ell$ and
$\algsl(\ell)_\aff$ at level $n$ implies the statements.
One of missed statements in \cite{Na-branching} is about the costalk
of $\IC$ of Uhlenbeck space. It will be deduced from analog of
\eqref{eq:1} in the affine setting.

In \cite{2016arXiv160403625B} Braverman-Finkelberg and the
second-named author generalized and refined the proposal in
\cite{braverman-2007} by replacing Uhlenbeck spaces by Coulomb
branches $\cM(\lambda,\mu)$ of quiver gauge theories.

Let us fix the notation for Coulomb branches $\cM(\lambda,\mu)$.  We
take a quiver $Q = (Q_0, Q_1)$, where $Q_0$ is the set of vertices,
and $Q_1$ is the set of oriented edges. Each oriented edge $h$ defines
the outgoing and incoming vertices, denoted by
$\vout{h}\xrightarrow{h}\vin{h}$.
We take two $Q_0$-graded vector spaces
\begin{equation*}
  V = \bigoplus_i V_i, \qquad W = \bigoplus_i W_i
\end{equation*}
and define $\lambda$, $\mu$ as $\lambda = \sum \dim W_i \Lambda_i$,
$\mu = \lambda - \sum \dim V_i\alpha_i$, where $\Lambda_i$ (resp.\
$\alpha_i$) is the $i$-th fundamental weight (resp.\ of simple root)
$\fg = \operatorname{Lie}G$. We then define
\begin{equation*}
  \mathbf G \defeq \prod_i \GL(V_i),\qquad
  \bN \defeq \bigoplus_h \Hom(V_{\vout{h}},V_{\vin{h}})
  \oplus\bigoplus_i \Hom(W_i,V_i).
\end{equation*}
We consider the variety $\cR$ of BFN triples for $\mathbf G$ (instead
of $G^\vee$), $\bN$ in the construction explained in \cref{sec:Coulomb}.
Then we define the Coulomb branch as
\begin{equation*}
  \cM(\lambda,\mu) \defeq \Spec H^{\mathbf G(\cO)}_*(\cR).
\end{equation*}

We have
\begin{Theorem}[\cite{2016arXiv160403625B}
  cf.~\cite{2019arXiv190706552N} for nonsymmetric case]
  Suppose $Q$ is of finite type and let $G^\vee$ be the corresponding
  group of adjoint type. Assume further that $\mu$ is dominant. Then
  $\cM(\lambda,\mu)\cong \oW^\lambda_{G^\vee,\mu}$.
\end{Theorem}

When $\mu$ is not a dominant, we need to introduce a variant of the
affine Grassmannian slice $\oW^\lambda_{G^\vee,\mu}$, called a
generalized slice. The above result holds for general $\mu$ if we replace
$\oW^\lambda_{G^\vee,\mu}$ by a generalized slice.

The generalized slice is expected to be isomorphic to the Uhlenbeck
partial compactification of the moduli space of singular
$G_c$-monopoles on $\RR^3$. Here we impose the usual asymptotic
behavior at $|x|\to\infty$, contrary to the result by
Braverman-Finkelberg \cite{braverman-2007} mentioned in
\cref{sec:monopole}.

Let us consider the case $G$ as above, dual of an untwisted affine
Kac-Moody group, and let $(G^\vee)_{\fin,c}$ as above. Then it is
expected that $\cM(\lambda,\mu)$ is the Uhlenbeck partial
compactification of $(G^\vee)_{\fin,c}$-instantons on the Taub-NUT
space divided by $\ZZ/\ell\ZZ$.
The Taub-NUT space is isomorphic to $\CC^2$ as a complex manifold, but
has a hyper-K\"ahler metric different from the Euclidean one.
The difference between instantons on $\RR^4$ and ones on the Taub-NUT
space corresponds to the difference of the asymptotic behavior of
monopoles at $|x|\to\infty$ mentioned above.

When $G$ is of affine type $A$, instantons on the Taub-NUT space are
described by variants of the ADHM description, which also involves a
certain nonlinear ordinary equation, called the Nahm's equation. This
is a result of Cherkis \cite{MR2525636}. He named instanton moduli
spaces \emph{bow varieties} in this description.

\begin{Theorem}[\cite{2016arXiv160602002N}]
  Suppose $Q$ is of affine type $A$. Then $\cM(\lambda,\mu)$ is
  isomorphic to a bow variety.
\end{Theorem}

Considering $S^1$-equivariant instantons on the Taub-NUT space, we can
get singular monopoles with the ordinary asymptotic behavior at
$|x|\to\infty$. Therefore the expectation above is true for arbitrary
$\mu$ when $G$ is of type $A$. Moreover when $\mu$ is dominant, a bow
variety is isomorphic to a quiver variety, which is the moduli space
of instantons on $\RR^4/(\ZZ/\ell\ZZ)$, as expected.

Let us return back to geometric Satake correspondence for a Kac-Moody
Lie group.

It is known that $\cM(\lambda,\mu)$ has an action of $T^\vee$.
It is conjectured that the fixed point set $\cM(\lambda,\mu)^{\rho=0}$
is either empty or a single point, and the latter happens if and only
if $V(\lambda)_\mu\neq 0$.
We consider the hyperbolic restriction functor $\Phi$ as above.
Let us denote by $\IC^\lambda_\mu$ the intersection cohomology complex
$\IC(\cM(\lambda,\mu))$ for brevity.
The proposal in \cite{2016arXiv160403625B} is an existence of a
natural isomorphism
\begin{equation}\label{eq:5}
  H^0(\Phi(\IC^\lambda_\mu)) \xrightarrow[\cong]{\text{?}} V(\lambda)_\mu.
\end{equation}

This generalization makes sense for a symmetric Kac-Moody Lie group
$G$ of simply-connected type, and can be even generalized to
symmetrizable non-symmetric cases thanks to
\cite{2019arXiv190706552N}. It is a refinement as Coulomb branches
make sense for any weights $\mu$, while Uhlenbeck spaces, as for
$\oW^\lambda_\mu$, correspond only to \emph{dominant} weights.
For finite type, this refinement is not so crucial, as we could use
$z^\mu\in\Gr_{G^\vee}$ for nondominant $\mu$. But it does not make
sense for general Kac-Moody Lie algebras.

The proposal in \cite{2016arXiv160403625B} was further refined by the
second-named author \cite{2018arXiv181004293N}: A concrete
construction of representations of $G$ was given under assumptions on
geometric properties of Coulomb branches. These assumptions were
proved for type $A^{(1)}_{n-1}$.

\section{Result -- in progress}\label{sec:main}

\subsection{}

Recall we have remarked a delicate difference between
$\Phi = \iota^* j^!$ and $i_\mu^! = \iota^! j^!$.
In order to talk a finer information, let us consider
$T^\vee$-equivariant cohomology groups as in \eqref{eq:2}.
There is a natural transformation $\iota^!\to\iota^*$, hence
\begin{equation}\label{eq:6}
  H^*_{T^\vee}(i_\mu^!\IC_{\mu}^\lambda)
  \to H^*_{T^\vee}(\Phi(\IC_{\mu}^\lambda)).
\end{equation}
The codomain is conjectured to be isomorphic to
$H^*_{T^\vee}(\mathrm{pt})\otimes
H^0(\Phi(\IC_{\mu}^\lambda))
\cong H^*_{T^\vee}(\mathrm{pt})\otimes V(\lambda)_\mu$,
where the second isomorphism is the conjectural geometric Satake.
The first isomorphism follows once we prove
$H^i(\Phi(\IC_{\mu}^\lambda))$ vanishes except in degree $i=0$ by a
standard argument. This vanishing follows from the conjectural
hyperbolic semi-smallness of the functor $\Phi$, which is a part of
conjectural statements to derive geometric Satake.

\eqref{eq:6} is a homomorphism between two free
$H^*_{T^\vee}(\mathrm{pt})$-modules of the same rank, but it is not an
isomorphism. 
By applying $(V(\lambda)\otimes\bullet)^G$ to the left hand side of
\eqref{eq:2} in the finite dimensional setting and using the Frobenius
reciprocity, we have
\begin{equation}
  \left(V(\lambda)\otimes \CC[(\fg/\fu)^*]\otimes \CC_{-\mu}\right)^B
\end{equation}
This makes sense in the Kac-Moody setting.

We consider a birational map $B/T\times\mathfrak t^*\to (\fg/\fu)^*$
given by the coadjoint action
$(bT,X) \mapsto \operatorname{Ad}_b^*(X)$. It is $B$-equivariant where
$B$ acts trivially on $\mathfrak t^*$.
We apply it to get get a homomorphism from the left hand side to
\begin{equation*}
 \left(V(\lambda)\otimes \CC[B/T\times \mathfrak t^*]\otimes
  \CC_{-\mu}\right)^B
= \CC[\mathfrak t^*]\otimes (V(\lambda)\times \CC_{-\mu})^T
= \CC[\mathfrak t^*]\otimes V(\lambda)_\mu. 
\end{equation*}
Recall that $\CC[\mathfrak{t}^*]$ is isomorphic to
$H^*_{T^\vee}(\mathrm{pt})$. Therefore the above is nothing but
$H^*_{T^\vee}(\mathrm{pt})\otimes V(\lambda)_\mu$, hence the codomain
of \eqref{eq:6}.
This view point was taken by Ginzburg-Riche \cite{MR3352528} first
in the finite dimensional setting. 

Let us state the main conjecture. It basically says that
Ginzburg-Riche's framework works in the Kac-Moody setting.
Namely there exists an isomorphism $\heartsuit$ of
$S(\mathfrak{t}) \cong H^*_{T^\vee}(\mathrm{pt})$-modules, which sits
in a commutative diagram
\begin{equation}\label{eq:7}
  \begin{tikzcd}
    \left(V(\lambda)\otimes \CC[(\fg/\fu)^*]\otimes \CC_{-\mu}\right)^B
    \arrow[r,"\heartsuit","\sim"']
    \arrow[d]
    & H^*_{T^\vee}(i_\mu^! \IC^\lambda_\mu)
    \arrow[d,"\eqref{eq:6}"]
    \\
    \CC[\mathfrak{t}^*]\otimes V(\lambda)_\mu
    \arrow[r,"\eqref{eq:5}","\sim"']
    & H^*_{T^\vee}(\Phi(\IC^\lambda_\mu))\rlap{.}
  \end{tikzcd}
\end{equation}

\subsection{Strategy of a proof}

Let us explain a strategy of a proof of the above conjecture. This
strategy was used in the finite dimensional case in \cite{MR3352528},
and will be also used in affine type $A$ in the work in progress with
D.~Muthiah.

Let $\CC[\ft_{\mathrm{rs}}^*]$ be the localization of $\CC[\ft^*]$
inverting all coroots. Thanks to the localization theorem for
equivariant cohomology groups, \eqref{eq:6} becomes an isomorphism
after $\otimes \CC[\ft_{\mathrm{rs}}^*]$.
On the other hand, the left vertical arrow in \eqref{eq:7} also
becomes an an isomorphism after $\otimes \CC[\ft_{\mathrm{rs}}^*]$.
Therefore $\heartsuit$ exists over $\otimes \CC[\ft_{\mathrm{rs}}^*]$.
Hence what we need to show is that $\heartsuit$ is an isomorphism over
$\CC[\ft^*]$. Possible poles of $\heartsuit$ appear at coroot
hyperplanes $\alpha^\vee = 0$.

Suppose $\alpha^\vee$ is a real coroot, i.e., a Weyl group translate
of a simple coroot $\alpha_i^\vee$. The Weyl group symmetry can be
realized in the diagram \eqref{eq:7} by the dynamical Weyl group
introduced by Etingof-Varchenko \cite{MR1901247}. This was shown in
\cite{MR3352528} for the finite dimensional setting, and expected to
be true in general.
Thanks to the dynamical Weyl group, $\heartsuit$ extends across
$\alpha^\vee = 0$ if and only if it extends across $\alpha^\vee_i=0$.
By the framework in \cite{2018arXiv181004293N} to define the geometric
Satake isomorphism \eqref{eq:5}, the behavior at $\alpha^\vee_i=0$
corresponds to the restriction to the Levi subgroup $L_i$
corresponding to $\alpha_i^\vee$. In particular, the extension of
$\heartsuit$ across $\alpha_i^\vee=0$ follows from the same assertion
for $L_i$. Then $L_i$ is basically $\SL_2$, one checks the extension
by an explicit calculation of all entries of \eqref{eq:7}.

If $G$ is a finite dimensional group, all coroots are real, hence
there is nothing more to be done.

If $G$ is of affine type, there is a primitive imaginary coroot $c
\begin{NB}
  = \sum a_i^\vee \alpha_i^\vee
\end{NB}%
$, and all imaginary roots are its multiple. Thus we need to
understand what happens at $c=0$.
Instead of $L_i$ above, we need to treat the (homogeneous) Heisenberg
subgroup $\Heis$.
\begin{NB}
  $\Heis$ is \verb+\Heis+. The notation should be changed.
\end{NB}%

On the other hand, in the right column of \eqref{eq:7}, we replace
$\cM(\lambda,\mu)$ by the fixed point set with respect to a
cocharacter $\rho_\delta\colon \CC^\times\to T^\vee$, whose
differential is contained in $c=0$, but not in other coroot
hyperplanes thanks to the localization theorem of equivariant
cohomology groups.
It is expected that the fixed point set
\(
  \cM(\lambda,\mu)^{\rho_\delta=0}
\)
is an appropriate symmetric power of $\cS_\ell = \CC^2/(\ZZ/\ell\ZZ)$.

In order to relate the symmetric power of $\cS_\ell$ with the
Heisenberg subgroup, we consider its resolution given by a Hilbert
scheme of points on the minimial resolution of $\cS_\ell$. Then one
uses the construction \cite{MR1441880}. Two entries of the right
column of~\eqref{eq:7} are easily understood by this construction, and
can be compared with the left column. This finihses the proof in
affine type $A$.

This strategy should work more general affine types. Most of
difficulties appear in the proof of geometric Satake \eqref{eq:5}. As
we remarked at the end of \cref{sec:KacMoody}, we need to check some
geometric properties of $\cM(\lambda,\mu)$.

\subsection*{Acknowledgment}

The research of HN was supported in part by the World Premier
International Research Center Initiative (WPI Initiative), MEXT,
Japan, and by JSPS Grant Number 16H06335.

\bibliographystyle{myamsalpha}
\bibliography{nakajima,mybib,coulomb}

\end{document}